\documentclass[a4paper, 11pt]{article}
\usepackage{my_preamble}

\title{\textbf{Complex algebraic stacks and morphisms of intersection complexes}}
\author{Matthew Huynh}
\date{July 30, 2024}

\begin{document}

\maketitle

\begin{abstract}
    We generalize a construction of Barthel-Brasselet-Fieseler-Gabber-Kaup in the setting of complex varieties to the setting of finite type, complex algebraic stacks. 
    Given two such stacks $\mc{X}, \mc{Y}$ with affine stabilizers, and a morphism between them, we construct a morphism from the pullback of the intersection complex of $\mc{Y}$ to the intersection complex of $\mc{X}$. 
    As an application, we show that the Borel-Moore fundamental class of a closed substack $\mc{Z}$ in a Deligne-Mumford stack $\mc{X}$ lifts to a class in the intersection cohomology of $\mc{X}$. 
\end{abstract}

\section{Introduction}
The rational cohomology groups of a smooth, complex projective variety enjoy many nice properties. They satisfy Poincar\'e duality, underlie pure Hodge structures, and form an $\mf{sl}_2(\C)$-representation via the Hard Lefschetz Theorem. Goresky and MacPherson defined intersection (co)homology groups in \cites{IH1, IH2} for a general class of topological spaces, and the intersection (co)homology groups of a singular projective variety enjoy the properties listed above. 

The original cycle-theoretic definition of intersection homology proceeds by specifying geometric cycles whose supports suitably intersect the singular strata of the space; these are called intersection cycles. On a complex variety $X$, it makes sense to ask: if we are given an algebraic cycle $Z$ in $X$, is there an intersection cycle $\gamma$ that is homologous to $Z$?

To answer this question, Barthel-Brasselet-Fieseler-Gabber-Kaup \cite{relevement} used the sheaf-theoretic description of intersection homology and the shifted intersection complex. Their main theorem is that for any morphism $f: X \to Y$ between two pure-dimensional complex varieties, there exists a morphism $\mu^f: f^{\ast}IC_Y \to IC_X$, (not necessarily unique), fitting into the commutative diagram
\begin{equation*}
    \xymatrix{
        f^{\ast}IC_Y \ar[r]^-{\mu^f} & IC_X \\
        f^{\ast}\Q_Y \ar[r]_{\sim} \ar[u]^{f^{\ast}c_Y} & \Q_X \ar[u]_{c_X},
    }
\end{equation*}
where $c_X: \Q_X \to IC_X$ is the natural comparison morphism. The morphisms they construct are called \textit{associated morphisms}. 

One global consequence of their main theorem for proper maps $f: X \to Y$ is the following. There exists a morphism of intersection homology groups, 
\[
    \nu_f: IH_{\ast}^{BM}(X) \to IH_{\ast}^{BM}(Y), 
\]
such that for any intersection homology class $[\gamma] \in IH_{\ast}^{BM}(X)$, the Borel-Moore homology class of any intersection cycle in $\nu_f([\gamma])$ is the Borel-Moore homology class of $f(\gamma)$. By applying this to the inclusion of an irreducible subvariety in a complex variety, they give a positive answer to the question posed earlier. Namely, if $Z$ is an irreducible subvariety of a pure-dimensional complex variety $X$, then there exists an intersection cycle $\gamma$ on $X$, whose Borel-Moore homology class coincides with the Borel-Moore fundamental class of $Z$. 

The purpose of this paper is to establish the existence of associated morphisms in the setting of complex algebraic stacks. We give a precise formulation of our main result, Theorem \ref{thm:main_theorem}, after introducing our notation and conventions. As an application, we show in \S 4 that for any closed, irreducible substack $\mc{Z}$ of a pure-dimensional, finite type, Deligne-Mumford stack $\mc{X}$ over $\C$, the Borel-Moore homology class of $\mc{Z}$ lifts to a class in the intersection cohomology of $\mc{X}$. 

\subsection{Notation and conventions}
\begin{enumerate}[label=(\roman*).]
    \item We use the convention of \cite{laumon1999champs} for algebraic stacks $\mc{X}$ over a scheme $S$. In other words, we assume that the diagonal $\Delta: \mc{X} \to \mc{X} \times_S \mc{X}$ is representable, separated, and quasi-compact. 
    
    From now on, a ``stack'' means an algebraic stack $\mc{X}$ over $\Spec(\C)$ of finite-type. 
    \item We use the categories and six-functor formalism for sheaves on stacks with adic coefficients developed in \cite{laszlo_olsson_2}. Given a stack $\mc{X}$, we denote by $\mathbf{D}^b_c(\mc{X},\Qlbar)$ the so-called bounded derived category of constructible complexes with $\Qlbar$-coefficients on $\mc{X}$. Given a morphism of stacks $f: \mc{X} \to \mc{Y}$, we denote by $f^{\ast},f_{\ast},f_!,f^!$ the derived functors between the categories $\mathbf{D}^b_c(\mc{X},\Qlbar)$ and $\mathbf{D}^b_c(\mc{Y},\Qlbar)$. 
    \begin{rmk}
        Suppose that $X$ is a separated scheme of finite type over $\Spec(\C)$. Then the category $\mathbf{D}^b_c(X,\Qlbar)$ is equivalent to the category defined by Deligne in \cite{weil2}*{I.I.2} (see \cite{laszlo_olsson_2}*{Theorem 3.1.6}). 
    \end{rmk}
    \item We use the theory of perverse sheaves on stacks developed in \cite{laszlo2009perverse} building on work in \cite{laszlo_olsson_2} and work only with the middle perversity. 
    If $\mc{X}$ is an irreducible stack, then $IC_{\mc{X}}$ denotes the intersection complex of $\mc{X}$, shifted by $-\dim\mc{X}$. For example, if $\mc{X}$ is also smooth, then $IC_{\mc{X}} \simeq (\Qlbar)_{\mc{X}}$. If $\mc{X}_1,\ldots,\mc{X}_r$ are the irreducible components of $\mc{X}$, and $\mu: \coprod_{i = 1}^r \mc{X}_i \to \mc{X}$ is the natural map, then $IC_{\mc{X}} := \mu_{\ast}(\oplus_i IC_{\mc{X}_i})$. 
    \item Given a stack $\mc{X}$, the natural comparison morphism $(\Qlbar)_{\mc{X}} \to IC_{\mc{X}}$ is denoted by $c_{\mc{X}}$. When $\mc{X}$ is smooth, this is an isomorphism. 
    \item Given a morphism $f:\mc{X} \to \mc{Y}$, the adjunction morphism $(\Qlbar)_{\mc{Y}} \to f_{\ast}(\Qlbar)_{\mc{X}}$ is denoted by $a_f$. 
\end{enumerate}

\subsection{Statement of the main theorem}
\begin{thm}
    \label{thm:main_theorem}
    (See \S 1.1 for our convention on stacks). Let $\mc{X}$ and $\mc{Y}$ be two stacks with affine stabilizers, and let $f:\mc{X} \to \mc{Y}$ be a morphism between them. Then there exists a morphism ${\mu^f: f^{\ast}IC_{\mc{Y}} \to IC_{\mc{X}}}$ making the following diagram commutative
    \begin{equation}
        \label{eq:assoc_morphism_diagram}
        \xymatrix{
            f^{\ast}IC_{\mc{Y}} \ar[r]^{\mu^f} & IC_{\mc{X}} \\
            f^{\ast}(\Qlbar)_{\mc{Y}} \ar[u]^{f^{\ast}c_{\mc{Y}}} \ar[r]_{\sim} & (\Qlbar)_{\mc{X}}. \ar[u]_{c_{\mc{X}}}
        }
    \end{equation}
\end{thm}
\begin{rmk}
    \label{rmk:adjoint}
    Since $(f^{\ast},f_{\ast})$ form an adjoint pair, the existence of a morphism $\mu^f: f^{\ast}IC_{\mc{Y}} \to IC_{\mc{X}}$ making diagram (\ref{eq:assoc_morphism_diagram}) commutative is equivalent to the existence of a morphism ${\lambda_f: IC_{\mc{Y}} \to f_{\ast}IC_{\mc{X}}}$ making the following diagram commutative
    \begin{equation}
        \label{eq:adjoint_assoc_morphism_diagram}
        \xymatrix{
            IC_{\mc{Y}} \ar[r]^{\lambda_f} & f_{\ast}IC_{\mc{X}} \\
            (\Qlbar)_{\mc{Y}} \ar[u]^{c_{\mc{Y}}} \ar[r]_{a_f} & f_{\ast}(\Qlbar)_{\mc{X}}. \ar[u]_{f_{\ast}c_{\mc{X}}}
        }
    \end{equation}
    Weber \cite{weber99} used this observation and the Decomposition Theorem for proper morphisms of complex varieties to give an alternative proof of \cite{relevement}*{Th\'eor\`eme 2.1}, and we adapt his argument to the setting of stacks. 

\end{rmk}
\begin{rmk}
    \label{rmk:uniqueness}
    The morphism $\mu^f$ is not necessarily unique. Given a morphism between two pure-dimensional complex varieties $f: X \to Y$, associated morphisms $\mu^f: f^{\ast}IC_{Y} \to IC_{X}$ are in bijection with lifts of the class $[\Gamma_f] \in H^{BM}_{2\dim X}(X\times Y, \Q)$ to $IH^{BM}_{2\dim X}(X \times Y, \Q)$, \cite{relevement}*{Th\'eor\`eme 2.7}. However, graphs of morphisms of stacks do not embed into the product of their source and target in general (consider the identity morphism $BG \to BG$ for a non-trivial group $G$. The graph of this morphism is the diagonal morphism $BG \to BG \times BG$, which is not a closed embedding). Therefore, it is not clear to the author if there is an analogous description of associated morphisms for morphisms of stacks. 
\end{rmk}
\subsection{Acknowledgements}
The author is thankful to 
Mark de Cataldo and to Andres Fernandez Herrero for helpful conversations and comments on earlier drafts of this paper, which have improved it greatly. The author is also thankful to the referee for helpful comments.

This paper was partially supported by N.S.F. Grant DMS-2200492 and by Simons Foundation International, LTD.

\section{A corollary of the Decomposition Theorem for stacks}
\begin{setup}
    \label{setup:IC_splits}
    Throughout this section, $\mc{A}$ and $\mc{B}$ denote stacks with affine stabilizers, and $\pi: \mc{A} \to \mc{B}$ is a proper, surjective morphism that is representable by schemes. 
\end{setup}

The purpose of this section is to prove an auxiliary result, which is a corollary of Sun's version of the Decomposition Theorem \cite{sun_DT_2016}*{Theorem 4.2.4}. The setting of Setup \ref{setup:IC_splits} satisfies the assumptions of the referenced theorem, and thus we have an isomorphism of objects in $\mathbf{D}^b_c(\mc{B},\Qlbar)$
\begin{equation}
    \label{eq:decomp_thm}
    \pi_{\ast}IC_{\mc{A}} \simeq \bigoplus_{i \in \Z} (\p{\ms{H}}^i\pi_{\ast}IC_{\mc{A}})[-i].
\end{equation}
Furthermore, the perverse sheaves $\p{\ms{H}}^i\pi_{\ast}IC_{\mc{A}}$ appearing in the right-hand side of (\ref{eq:decomp_thm}) are semisimple, and hence split into direct sums of simple perverse sheaves. We can deduce the following
\begin{prop}
    \label{prop:cor_of_DT}
    Assume Setup \ref{setup:IC_splits}. Then the complex $IC_{\mc{B}}$ is a direct summand of $\pi_{\ast}IC_{\mc{A}}$. 
\end{prop}
\begin{rmk}
    A proof of this proposition for morphisms of varieties over finite fields, (which can be adapted to the setting of complex varieties) can be found in \cite{de2018frobenius}*{\S 5.3}, and our proof strategy follows theirs. The proposition does not seem to immediately follow from the case of schemes and descent, because the constructed splitting is a morphism of objects in the derived category. 
\end{rmk}

\begin{proof}[\proofname\ of Proposition \ref{prop:cor_of_DT}]
By working with irreducible components one at a time, we may assume that $\mc{A}$ and $\mc{B}$ are irreducible. Furthermore, $\mc{A}$ (resp. $\mc{B}$) and $\mc{A}_{\rm{red}}$ (resp. $\mc{B}_{\rm{red}}$) have equivalent lisse-\'etale sites, so we may assume that $\mc{A}$ and $\mc{B}$ are reduced. 
First, we record two lemmas, which let us deduce the proposition from a special case.

\begin{lem}
    \label{lem:reduce_to_open}
    Let $j: \mc{B}^{\circ} \to \mc{B}$ be the inclusion of an open, dense substack, and denote by $\pi^{\circ}: \mc{A}^{\circ} \to \mc{B}^{\circ}$ the pullback of $\pi$ to $\mc{B}^{\circ}$. Assume that $IC_{\mc{B}^{\circ}}$ is a direct summand of $\pi^{\circ}_{\ast}IC_{\mc{A}^{\circ}} = j^{\ast}\pi_{\ast}IC_{\mc{A}}$. Then $IC_{\mc{B}}$ is a direct summand of $\pi_{\ast}IC_{\mc{A}}$. 
    \begin{proof}
        Recall that every simple perverse sheaf on $\mc{B}$ is obtained by: (i) taking an irreducible smooth $\Qlbar$-sheaf $L$ on a smooth, locally closed, and irreducible substack $\mc{V}$ of $\mc{B}$; (ii) shifting $L$ by $\dim \mc{V}$; (iii) applying the intermediate extension functor to obtain a perverse sheaf on $\ol{\mc{V}}$; (iv) pushing forward to $\mc{B}$ (see \cite{laszlo2009perverse}*{Theorem 8.2}). 
        
        Moreover, if we have an open, dense substack $\mc{U}$ inside $\mc{V}$, then the intermediate extension of $(L[\dim\mc{V}])|_{\mc{U}}$ to $\ol{\mc{V}}$ followed by pushing forward to $\mc{B}$ yields the same object as the one described above. This follows straightforwardly from properties of the intermediate extension functor (see \cite{laszlo2009perverse}*{Lemma 6.1}). 

        The shifted simple summands of $\pi_{\ast}IC_{\mc{A}}$ (resp. $\pi^{\circ}_{\ast}IC_{\mc{A}^{\circ}}$) are unique. Thus by the assumption of the lemma, there is a simple perverse sheaf on $\mc{B}$ shifted by $-\dim\mc{B}^{\circ}$, which is a summand of $\pi_{\ast}IC_{\mc{A}}$, and whose restriction to $\mc{B}^{\circ}$ is isomorphic to $IC_{\mc{B}^{\circ}}$. Call this object $IC_{\ol{\mc{V}}}(L)$. It follows that $\ol{\mc{V}} = \mc{B}$, since $IC_{\mc{B}^{\circ}}$ is supported on all of $\mc{B}^{\circ}$. 

        Since $IC_{\mc{B}}(L)$ restricted to $\mc{V} \cap \mc{B}^{\circ}$ is $L|_{\mc{V}\cap\mc{B}^{\circ}}$, and is isomorphic to $IC_{\mc{B}^{\circ}}$ restricted to $\mc{V}\cap \mc{B}^{\circ}$, it follows that $L|_{\mc{V}\cap\mc{B}^{\circ}} \simeq (\Qlbar)_{\mc{V}\cap\mc{B}^{\circ}}$. Hence $IC_{\mc{B}}(L) = IC_{\mc{B}}$.  
    \end{proof}
\end{lem}

\begin{lem}
    \label{lem:reduce_to_normal}
    Let $\nu: \mc{A}^{\nu} \to \mc{A}$ be the normalization of $\mc{A}$ (see \cite{stacks-project}*{\href{https://stacks.math.columbia.edu/tag/0GMH}{Tag 0GMH}}). If $IC_{\mc{B}}$ is a direct summand of $(\pi \circ \nu)_{\ast}IC_{\mc{A}^{\nu}}$, then it is also a direct summand of $\pi_{\ast}IC_{\mc{A}}$. 
    \begin{proof}
        It suffices to show $\nu_{\ast}IC_{\mc{A}^{\nu}} \simeq IC_{\mc{A}}$, from which the lemma follows immediately. Let $j: \mc{U} \hookrightarrow \mc{A}$ be the inclusion of an open, dense, smooth substack over which the morphism $\nu$ is an isomorphism, and let $i: \mc{Z} \hookrightarrow \mc{A}$ be the inclusion of the closed complement of $\mc{U}$. The situation fits into a commutative diagram
        \[
            \xymatrix{
                \mc{Z} \times_{\mc{A}} \mc{A}^{\nu} \ar[r]^-{i'} \ar[d]_{\nu_{\mc{Z}}} & \mc{A}^{\nu} \ar[d]_{\nu} & \mc{U}^{\nu} \ar[d]^{\nu_{\mc{U}}}_{\sim} \ar[l] \\
                \mc{Z} \ar[r]_{i} & \mc{A} & \mc{U}. \ar[l]^{j}
            }
        \]
        Because $IC_{\mc{A}}[\dim\mc{A}]$ is the intermediate extension of $(\Qlbar)_{\mc{U}}[\dim\mc{U}]$, it suffices to show that 
        \[
            \xymatrix{
                j^{\ast}\nu_{\ast}IC_{\mc{A}^{\nu}}[\dim\mc{A}] \simeq (\Qlbar)_{\mc{U}}[\dim\mc{U}]   
            }\hspace{20pt}\text{and}\hspace{20pt}
            \xymatrix{
                \p{\ms{H}}^0i^{\ast}(\nu_{\ast}IC_{\mc{A}^{\nu}}[\dim\mc{A}])= 0,
            }
        \]
        \cite{laszlo2009perverse}*{Lemma 6.1}. The first equality holds because $IC_{\mc{A}^{\nu}}[\dim\mc{A}]$ is the intermediate extension of $(\Qlbar)_{\mc{U}^{\nu}}[\dim\mc{U}^{\nu}]$, and by smooth base change, 
        \[ 
            j^{\ast}\nu_{\ast}IC_{\mc{A}^{\nu}}[\dim\mc{A}] \simeq (\nu_{\mc{U}})_{\ast}(\Qlbar)_{\mc{U}^{\nu}}[\dim\mc{U}^{\nu}] \simeq (\Qlbar)_{\mc{U}}[\dim\mc{U}].
        \]
        As for the second equality, the pushforward by a finite, schematic morphism of stacks is perverse $t$-exact, which follows from the analogous statement for varieties, \cite{bbd}*{Corollaire 4.1.3}, and the definition of the perverse $t$-structure on a stack, \cite{laszlo2009perverse}*{\S 4}. Therefore, we have
        \[
            \p{\ms{H}}^0i^{\ast}(\nu_{\ast}IC_{\mc{A}^{\nu}}[\dim\mc{A}]) = \p{\ms{H}}^0(\nu_{\mc{Z}})_{\ast}(i')^{\ast}IC_{\mc{A}^{\nu}}[\dim\mc{A}] = (\nu_{\mc{Z}})_{\ast}\p{\ms{H}}^0(i')^{\ast}IC_{\mc{A}^{\nu}} = 0.
        \]
    \end{proof}
\end{lem}

\begin{rmk}
    We would like to thank the anonymous referee for bringing the following historical remark to our attention.

    In the development of intersection cohomology, Goresky and MacPherson were seeking a definition of cohomology groups that would not change after normalization. Therefore, the isomorphism $\nu_{\ast}IC_{A^{\nu}} \simeq IC_A$, (where $\nu: A^{\nu} \to A$ is the normalization of an integral complex variety), was an expected property of intersection cohomology. 
\end{rmk}

    By replacing $\mc{B}$ with an open, dense substack and applying Lemma \ref{lem:reduce_to_open}, we may assume that each simple summand appearing in the decomposition of $\pi_{\ast}IC_{\mc{A}}$ has full support, each shifted simple summand is actually a shifted irreducible smooth $\Qlbar$-sheaf, and $\mc{B}$ is smooth. Note that in this case, $\pi_{\ast}IC_{\mc{A}} \simeq \oplus_i (R^i\pi_{\ast}IC_{\mc{A}})[-i]$. This reduces the proposition to showing that $(\Qlbar)_{\mc{B}}$ is a direct summand of $R^0\pi_{\ast}IC_{\mc{A}}$. 

    By Lemma \ref{lem:reduce_to_normal}, we may assume that $\mc{A}$ is normal, which implies that $R^0\pi_{\ast}IC_{\mc{A}} \simeq R^0\pi_{\ast}(\Qlbar)_{\mc{A}}$. To see this, let $U \to \mc{A}$ be a smooth presentation with $U$ normal. By smooth base change, the morphism $\ms{H}^0(c_{\mc{A}}): (\Qlbar)_{\mc{A}} \to \ms{H}^0IC_{\mc{A}}$ pulls back to the isomorphism $\ms{H}^0(c_U): (\Qlbar)_U \to \ms{H}^0IC_{U}$. Furthermore, since the morphism $c_U$ is compatible with smooth base change, this isomorphism descends to an isomorphism $(\Qlbar)_{\mc{A}} \to \ms{H}^0IC_{\mc{A}}$. This yields a distinguished triangle $((\Qlbar)_{\mc{A}},IC_{\mc{A}},\tau_{\geq 1}IC_{\mc{A}})$, and pushing forward by $\pi$ yields the claim. 

    It only remains to show that $(\Qlbar)_{\mc{B}}$ is a direct summand of $R^0\pi_{\ast}(\Qlbar)_{\mc{A}}$, after possibly replacing $\mc{B}$ by an open, dense substack. Let $p: U \to \mc{B}$ be a smooth presentation, and let $\pi_U: \mc{A}_U \to U$ denote the base change. The map $\pi_U$ is a proper morphism of schemes, and hence has a Stein factorization
    \[
        \xymatrix{
            \mc{A}_U \ar[r]^{\pi_U} \ar[dr]_{g_U} & U \\
            & U', \ar[u]_{h_U}
        }
    \]
    where $g_U$ is proper with connected fibers, and $h_U$ is finite. Next, let $V$ be an open, dense subscheme of $U$ such that over $V$, the finite morphism $h_U$ is unramified. Replace $U$ by $V$, and replace $\mc{B}$ by the image of $V$ in $\mc{B}$.

    Observe that $p^{\ast}(\Qlbar)_{\mc{B}}$ is canonically isomorphic to $(\Qlbar)_U$, and $p^{\ast}R^0\pi_{\ast}(\Qlbar)_{\mc{A}}$ is isomorphic to $R^0(\pi_U)_{\ast}(\Qlbar)_{\mc{A}_U}$ via the base change morphism, since $p$ is smooth. Therefore, the pullback of $\ms{H}^0(a_{\pi}): (\Qlbar)_{\mc{B}} \to R^0\pi_{\ast}(\Qlbar)_{\mc{A}}$ by $p$ becomes $\ms{H}^0(a_{\pi_U}): (\Qlbar)_U \to R^0(\pi_U)_{\ast}(\Qlbar)_{\mc{A}_U}$, after making these identifications. Conversely, the morphism $\ms{H}^0(a_{\pi_U})$ descends to $\mc{B}$. 

    Define $\varphi$ to be the composition
    \[
        \xymatrix{
            R^0(\pi_U)_{\ast}(\Qlbar)_{\mc{A}_U} = R^0(h_U)_{\ast}R^0(g_U)_{\ast}(\Qlbar)_{\mc{A}_U} \ar[d]^{R^0(h_U)_{\ast}\ms{H}^0(a_{g_U})^{-1}} \\ R^0(h_U)_{\ast}(\Qlbar)_{U'} \ar[d]^{(1/\mathrm{deg}(h_U))\tr(h_U)} \\ (\Qlbar)_U,
        }
    \]
    where $\ms{H}^0(a_{g_U}): (\Qlbar)_U \to R^0(g_U)_{\ast}(\Qlbar)_{U'}$ is an isomorphism, since $g_U$ has connected fibers, and $\tr(h_U): (h_U)_{\ast}h_U^{\ast}(\Qlbar)_U \to (\Qlbar)_U$ is the trace morphism defined in \cite{sga4vol3}*{Th\'eor\`eme 6.2.3, Expos\'e XVII}. Then the composition $\varphi \circ \ms{H}^0(a_{\pi_U})$ is the identity morphism. 
    
    It only remains to show that $\varphi$ descends to a morphism of sheaves on $\mc{B}$. To this end, let $s,t: R:= U \times_{\mc{B}} U \rightrightarrows U$ be the groupoid of relations. Let $\mc{A}_R$ denote 
    \[
        \mc{A}_U \times_{\pi_U,U,s} R \cong \mc{A} \times_{\pi,\mc{B},s\circ\pi} R \cong \mc{A} \times_{\pi,\mc{B},t\circ \pi}R \cong \mc{A}_U \times_{\pi_U,U,t} R, 
    \]
    and let $\pi_R: \mc{A}_R \to R$ denote the base change. Both factors of $\varphi$ are compatible with smooth base change, so the pullback of $\varphi$ by $s$ and $t$ is identified with the composition \[
        \xymatrix{
            R^0(\pi_R)_{\ast}(\Qlbar)_{\mc{A}_R} = R^0(h_R)_{\ast}R^0(g_R)_{\ast}(\Qlbar)_{\mc{A}_R} \ar[d]^{R^0(h_R)_{\ast}\ms{H}^0(a_{g_R})^{-1}} \\
            R^0(h_R)_{\ast}(\Qlbar)_{R'} \ar[d]^{(1/\mathrm{deg}(h_R))\tr(h_R)} \\
            (\Qlbar)_R,
        }
    \]
    where 
    \[ 
        \xymatrix{
            \mc{A}_R \ar[r]^-{g_R} & R' \cong U' \times_{g_U,U,s} R \cong U'\times_{g_U,U,t} R\ar[r]^-{h_R} & R
        }
    \]
    is the Stein factorization of $\pi_R: \mc{A}_R \to R$. Therefore, $\varphi$ descends to a morphism of sheaves on $\mc{B}$. 
\end{proof}

\begin{rmk}
    \label{rmk:special_case}
    \, 
    
    \begin{enumerate}[label=(\roman*).]
        \item Suppose in addition to Setup \ref{setup:IC_splits} that $\mc{A}$ and $\mc{B}$ are smooth and irreducible. Then the injection $IC_{\mc{B}} \to \pi_{\ast}IC_{\mc{A}}$ constructed in the proof of Proposition \ref{prop:cor_of_DT} has a simple description. It is the composition 
        \[
            \xymatrixcolsep{3pc}\xymatrix{
                IC_{\mc{B}} \ar[r]^-{c_{\mc{B}}^{-1}}_-{\sim} & (\Qlbar)_{\mc{B}} \ar[r]^-{a_{\pi}} & \pi_{\ast}(\Qlbar)_{\mc{A}} \ar[r]^-{\pi_{\ast}c_{\mc{A}}}_-{\sim} & \pi_{\ast}IC_{\mc{A}}.
            }
        \]
        \item The splitting 
        \begin{equation}
            \label{eq:splitting}
            IC_{\mc{B}} \to \pi_{\ast}IC_{\mc{A}} \to IC_{\mc{B}},
        \end{equation}
        constructed in the proof of Proposition \ref{prop:cor_of_DT} is compatible with restriction in the following sense. If $\mc{B}$ is irreducible, $j: \mc{B}^{\circ} \hookrightarrow \mc{B}$ is the inclusion of an open substack, and $\pi^{\circ}: \mc{A}^{\circ} \to \mc{B}^{\circ}$ is the base change of $\pi$, then the splitting
        \[
            IC_{\mc{B}^{\circ}} \to \pi^{\circ}_{\ast}IC_{\mc{A}^{\circ}} \to IC_{\mc{B}^{\circ}},
        \]
        provided by Proposition \ref{prop:cor_of_DT} \textit{is} the pullback of (\ref{eq:splitting}) by $j$, after we make the identifications $IC_{\mc{B}^{\circ}} \simeq j^{\ast}IC_{\mc{B}}$ and $\pi^{\circ}_{\ast}IC_{\mc{A}^{\circ}} \simeq j^{\ast}\pi_{\ast}IC_{\mc{A}}$ via smooth base change. 
    \end{enumerate}
\end{rmk}

\section{Proof of Theorem \ref{thm:main_theorem}}
Recall the setting of the main theorem, i.e. $\mc{X}, \mc{Y}$ are two stacks with affine stabilizers, and $f: \mc{X} \to \mc{Y}$ is a morphism between them. Following Remark \ref{rmk:adjoint}, we will show there exists a morphism $\lambda_f: IC_{\mc{Y}} \to f_{\ast}IC_{\mc{X}}$ such that $\lambda_f$ makes diagram (\ref{eq:adjoint_assoc_morphism_diagram}) commutative.
\begin{proof}[\proofname\ of Theorem \ref{thm:main_theorem}]
    We may assume that $\mc{X}$ and $\mc{Y}$ are irreducible and reduced, by the same reasoning provided at the beginning of the proof of Proposition \ref{prop:cor_of_DT}. 

Let $\pi: \widetilde{\mc{Y}} \to \mc{Y}$ be the resolution of singularities constructed in \cite{resolution_stacks}*{Theorem 5.1.1}. Then the morphism $\pi$ is proper and representable by schemes, since the resolution is obtained by blowing up closed substacks. In particular, $\widetilde{\mc{Y}}$ has affine stabilizers. 

Consider the cartesian diagram
\begin{equation*}
    \xymatrix{
        \widetilde{\mc{X}} \ar[r]^{\tilde{f}} \ar[d]_{\pi'} & \widetilde{\mc{Y}} \ar[d]^{\pi} \\
        \mc{X} \ar[r]_{f} & \mc{Y}.\ar@{}[ul]|\square
    }
\end{equation*}
Then $\pi'$ is also proper and representable by schemes, and $\widetilde{X}$ has affine stabilizers. By applying the Decomposition Theorem and Proposition \ref{prop:cor_of_DT} to $\pi$ and $\pi'$, we obtain morphisms $i: IC_{\mc{Y}} \to \pi_{\ast}IC_{\widetilde{\mc{Y}}}$ and $p: \pi'_{\ast}IC_{\widetilde{\mc{X}}} \to IC_{\mc{X}}$ fitting into a diagram
\begin{equation}
    \label{eq:big_diagram}
    \xymatrixcolsep{3.5pc}\xymatrix{
        \pi_{\ast}IC_{\widetilde{\mc{Y}}} \ar[r]^{\pi_{\ast}(c_{\widetilde{\mc{Y}}})^{-1}} & \pi_{\ast}(\Qlbar)_{\widetilde{\mc{Y}}} \ar[r]^-{\pi_{\ast}a_{f'}} & f_{\ast}\pi'_{\ast}(\Qlbar)_{\widetilde{\mc{X}}} \ar[r]^-{f_{\ast}\pi'_{\ast}c_{\widetilde{\mc{X}}}} & f_{\ast}\pi'_{\ast}IC_{\widetilde{\mc{X}}} \ar[d]_{f_{\ast}p} \\
        IC_{\mc{Y}} \ar[u]^{i} & (\Qlbar)_{\mc{Y}} \ar[l]^{c_{\mc{Y}}} \ar[r]_{a_f} \ar[u]^{a_{\pi}} & f_{\ast}(\Qlbar)_{\mc{X}} \ar[r]_{f_{\ast}c_{\mc{X}}} \ar[u]^{f_{\ast}a_{\pi'}} & f_{\ast}IC_{\mc{X}}.
    }
\end{equation}
The desired morphism $\lambda_f$ is the composition
\[
    \xymatrixcolsep{2.4pc}\xymatrix{
        IC_{\mc{Y}} \ar[r]^{i} & \pi_{\ast}IC_{\tilde{\mc{Y}}} \ar[r]^{\pi_{\ast}(c_{\tilde{\mc{Y}}})^{-1}} & \pi_{\ast}(\Qlbar)_{\tilde{\mc{Y}}} \ar[r]^{\pi_{\ast}a_{f'}} & f_{\ast}\pi_{\ast}'(\Qlbar)_{\tilde{\mc{X}}} \ar[r]^{f_{\ast}\pi_{\ast}'c_{\tilde{\mc{X}}}} & f_{\ast}\pi_{\ast}'IC_{\tilde{\mc{X}}} \ar[r]^{f_{\ast}p} & f_{\ast}IC_{\mc{X}},
    }
\]
and we must show diagram (\ref{eq:big_diagram}) is commutative. 

The middle square is clearly commutative, so it suffices to establish the commutativity of the two squares below
\begin{equation}
    \label{eq:two_squares}
    \xymatrix{
        \pi_{\ast}IC_{\widetilde{\mc{Y}}} \ar[r]^{\pi_{\ast}(c_{\widetilde{\mc{Y}}})^{-1}} & \pi_{\ast}(\Qlbar)_{\widetilde{\mc{Y}}} \\
        IC_{\mc{Y}} \ar[u]^{i} & (\Qlbar)_{\mc{Y}} \ar[l]^{c_{\mc{Y}}} \ar[u]_{a_{\pi}}
    }
    \hspace{40pt}%
    \xymatrix{
        \pi'_{\ast}(\Qlbar)_{\widetilde{\mc{X}}} \ar[r]^-{\pi'_{\ast}c_{\widetilde{\mc{X}}}} & \pi'_{\ast}IC_{\widetilde{\mc{X}}} \ar[d]^{p} \\
        (\Qlbar)_{\mc{X}} \ar[r]_{c_{\mc{X}}} \ar[u]^{a_{\pi'}} & IC_{\mc{X}}.
    }
\end{equation}

Let $j_{\mc{U}}:\mc{U} \hookrightarrow \mc{Y}$ (resp. $j_{\mc{V}}: \mc{V} \hookrightarrow \mc{X}$) be the inclusion of an open, dense, smooth substack. After we pullback the left square of (\ref{eq:two_squares}) by $j_{\mc{U}}$, it commutes by part (i) of Remark \ref{rmk:special_case}. Similarly, after we pullback the right square of (\ref{eq:two_squares}) by $j_{\mc{V}}$, it commutes by part (ii) of Remark \ref{rmk:special_case}.

It remains to show that if two morphisms from $(\Qlbar)_{\mc{Y}}$ to $\pi_{\ast}(\Qlbar)_{\widetilde{\mc{Y}}}$ agree after restriction to $\mc{U}$, then they must coincide before restriction (resp. if two morphisms from $(\Qlbar)_{\mc{X}}$ to $IC_{\mc{X}}$ agree after restriction to $\mc{V}$, then they must coincide before restriction). To this end, observe that the restriction functor $j_{\mc{U}}^{\ast}$ yields an isomorphism
\[
    \xymatrix{
        \Hom_{\mathbf{D}^b_c(\mc{Y},\Qlbar)}((\Qlbar)_{\mc{Y}},\pi_{\ast}(\Qlbar)_{\widetilde{\mc{Y}}}) = H^0(\widetilde{\mc{Y}},(\Qlbar)_{\widetilde{\mc{Y}}}) \ar[d]^{j_{\mc{U}}^{\ast}} \\
        \Hom_{\mathbf{D}^b_c(\mc{U},\Qlbar)}((\Qlbar)_{\mc{U}},\pi_{\ast}(\Qlbar)_{\widetilde{\mc{U}}}) = H^0(\widetilde{\mc{U}},(\Qlbar)_{\widetilde{\mc{U}}}),
    }
\]
because it is simply the restriction of global sections of the constant sheaf. Similarly, the restriction functor $j_{\mc{V}}^{\ast}$ yields an isomorphism
\[
    \xymatrix{
        \Hom_{\mathbf{D}^b_c(\mc{X},\Qlbar)}((\Qlbar)_{\mc{X}},IC_{\mc{X}}) = H^0(\mc{X},IC_{\mc{X}}) \cong \Qlbar \ar[d]^{j_{\mc{V}}^{\ast}} \\
        \Hom_{\mathbf{D}^b_c(\mc{V},\Qlbar)}((\Qlbar)_{\mc{V}},IC_{\mc{V}}) = H^0(\mc{V},IC_{\mc{V}}) \cong \Qlbar,
    }
\]
since it sends the comparison morphism $c_{\mc{X}}$ to $c_{\mc{V}}$. Thus the theorem is proven.
\end{proof}

\section{Lifts of classes of closed substacks}
\begin{setup}
    \label{setup:lift_of_class}
    Throughout this section, $\mc{X}$ denotes a Deligne-Mumford stack with affine stabilizers, which has pure-dimension $n$. $\mc{Z}$ denotes an $m$-dimensional irreducible, closed substack of $\mc{X}$, and the inclusion is denoted $i: \mc{Z} \hookrightarrow \mc{X}$.
\end{setup}
Here we present an application of Theorem \ref{thm:main_theorem}, which is a generalization of \cite{relevement}*{Th\'eor\`eme 2.4}. Recall that the Borel-Moore homology groups of $\mc{X}$ are defined by 
\[
    H^{BM}_{q}(\mc{X}) := H^{-q}(\mc{X},\omega_{\mc{X}}),
\]
where $\omega_{\mc{X}}$ is the dualizing complex constructed in \cite{laszlo_olsson_2}. These groups are zero when $q \notin [0,2n]$, and the rank of $H^{BM}_{2n}(\mc{X})$ equals the number of irreducible components of $\mc{X}$. In particular, $H^{BM}_{2m}(\mc{Z})$ is generated by a single element $[\mc{Z}]$, called the fundamental class. For more on the basic properties of Borel-Moore homology for Deligne-Mumford stacks, see \cite{olsson_BMhomology}*{\S 2}. 

The goal of this section is to show that the image of $[\mc{Z}]$ in $H^{BM}_{2m}(\mc{X})$ lifts to a class in the intersection cohomology of $\mc{X}$. By applying Verdier duality to the comparison morphism $c_{\mc{Z}}: (\Qlbar)_{\mc{Z}} \to IC_{\mc{Z}}$ and then shifting by $-2m$, we obtain a morphism $IC_{\mc{Z}} \to \omega_{\mc{Z}}[-2m]$. This induces an isomorphism on cohomology in degree 0, as we show below. 

\begin{lem}
    \label{lem:cohomology_iso}
    Assume Setup \ref{setup:lift_of_class}. Then the morphism $IC_{\mc{Z}} \to \omega_{\mc{Z}}[-2m]$, which is obtained by dualizing the comparison morphism and shifting, induces an isomorphism 
    \[
        H^0(\mc{Z},IC_{\mc{Z}}) \isomto H^{-2m}(\mc{Z},\omega_{\mc{Z}}).
    \]
    \begin{proof}
        Let $j: \mc{U} \hookrightarrow \mc{Z}$ be the inclusion of an open, dense, smooth substack, and let $i: \mc{W} \hookrightarrow \mc{Z}$ be the inclusion of the closed complement. Apply the functorial distinguished triangle $i_!i^! \to \id \to j_{\ast}j^{\ast}$ to the composition $(\Qlbar)_{\mc{Z}} \to IC_{\mc{Z}} \to \omega_{\mc{Z}}[-2m]$ and then apply the functor $H^0(\mc{Z},-)$. This yields a commutative diagram
        \[
            \xymatrixcolsep{1.5pc}\xymatrix{
                H^0(\mc{W},i^!(\Qlbar)_{\mc{Z}}) \ar[r] \ar[d] & H^0(\mc{Z},(\Qlbar)_{\mc{Z}}) \ar[r]^-{\simeq} \ar[d] & H^0(\mc{Z},j_{\ast}(\Qlbar)_{\mc{U}}) \ar[r] \ar[d]_{\simeq} & H^1(\mc{W},i^!(\Qlbar)_{\mc{Z}}) \ar[d] \\
                H^0(\mc{W},i^!IC_{\mc{Z}}) \ar[r] \ar[d] & H^0(\mc{Z},IC_{\mc{Z}}) \ar[r] \ar[d] & H^0(\mc{Z},j_{\ast}(\Qlbar)_{\mc{U}}) \ar[r] \ar[d]_{\simeq} & H^1(\mc{W},i^!IC_{\mc{Z}}) \ar[d] \\
                H^{-2m}(\mc{W},\omega_{\mc{W}}) \ar[r] & H^{-2m}(\mc{Z},\omega_{\mc{Z}}) \ar[r] & H^0(\mc{Z},j_{\ast}(\Qlbar)_{\mc{U}}) \ar[r] & H^{-2m+1}(\mc{W},\omega_{\mc{W}}),
            }
        \]
        where the rows are exact. Because $\dim\mc{W} \leq m-1$, we have 
        \[ 
            H^{-2m}(\mc{W},\omega_{\mc{W}}) = H^{-2m+1}(\mc{W},\omega_{\mc{W}}) = 0.
        \]
        By exactness, the morphism $H^{-2m}(\mc{Z},\omega_{\mc{Z}}) \to H^0(\mc{Z},j_{\ast}(\Qlbar)_{\mc{U}})$ is an isomorphism, so it suffices to show the surjective morphism $H^0(\mc{Z},IC_{\mc{Z}}) \to H^0(\mc{Z},j_{\ast}(\Qlbar)_{\mc{U}})$ is also injective. 

        Observe that by \cite{bbd}*{Corollaire 1.4.24} and \cite{laszlo2009perverse}*{Lemma 6.2}, we have $i^!IC_{\mc{Z}} \simeq \p\tau_{\geq m+1}i^!IC_{\mc{Z}}$, so $\p{\ms{H}}^s(i^!IC_{\mc{Z}}) = 0$ for all $s \leq m$. Thus the perverse spectral sequence 
        \[
            E_1^{r,s} = H^r(\mc{W},\p{\ms{H}}^s(i^!IC_{\mc{Z}})) \implies H^{r+s}(\mc{W},i^!IC_{\mc{Z}}),
        \]
        shows that $H^0(\mc{W},i^!IC_{\mc{Z}}) = 0$, and the lemma follows.
    \end{proof}
\end{lem}

\begin{rmk}
    Note that the proof of Lemma \ref{lem:cohomology_iso} also shows that the canonical map $H^0(\mc{Z},(\Qlbar)_{\mc{Z}}) \to H^0(\mc{Z},IC_{\mc{Z}})$ is an isomorphism. 
\end{rmk}

By applying the main theorem to the inclusion $i: \mc{Z} \to \mc{X}$, we can deduce the following
\begin{thm}
    \label{thm:lifts}
    Assume Setup \ref{setup:lift_of_class}. Then there is a class in $H^{2(n-m)}(\mc{X},IC_{\mc{X}})$ that maps to the image of $[\mc{Z}]$ in $H^{BM}_{2m}(\mc{X})$. 
    \begin{proof}
        By Theorem \ref{thm:main_theorem}, there exists an associated morphism $\mu^i: i^{\ast}IC_{\mc{X}} \to IC_{\mc{Z}}$. After applying Verdier duality on $\mc{Z}$ and shifting, we obtain the following commutative diagram
        \begin{equation}
            \label{eq:fund_class_diagram1}
            \xymatrix{
                IC_{\mc{Z}} \ar[r]^-{\mathbb{D}(\mu^i)} \ar[d] & i^!IC_{\mc{X}}[2(n-m)] \ar[d] \\
                \omega_{\mc{Z}}[-2m] \ar[r] & i^!\omega_{\mc{X}}[-2m].
            }
        \end{equation}
       Recall that there is a natural transformation $H^0(\mc{Z},i^!(-)) \to H^0(\mc{X},-)$, which comes from the adjunction morphism $i_!i^! \to \id$. Thus from (\ref{eq:fund_class_diagram1}) we obtain the commutative diagram
       \begin{equation}
           \label{eq:fund_class_diagram2}
           \xymatrix{
                H^0(\mc{Z},IC_{\mc{Z}}) \ar[r] \ar[d] & H^{2(n-m)}(\mc{Z},i^!IC_{\mc{X}}) \ar[r] \ar[d] & H^{2(n-m)}(\mc{X},IC_{\mc{X}}) \ar[d] \\
                H^{-2m}(\mc{Z},\omega_{\mc{Z}}) \ar[r]_-{=} & H^{-2m}(\mc{Z},i^!\omega_{\mc{X}}) \ar[r] & H^{-2m}(\mc{X},\omega_{\mc{X}}).
           }
       \end{equation}
       The left vertical arrow is an isomorphism by Lemma \ref{lem:cohomology_iso}, and the theorem follows.
    \end{proof}
\end{thm}

\bibliography{references}
{\footnotesize
    \textsc{Department of Mathematics, Stony Brook University,
        Stony Brook, NY 11794-3651,
    USA}\par\nopagebreak
      \textit{E-mail address}: \texttt{matthew.huynh@stonybrook.edu}
}

\end{document}